\documentclass{article}
\newcommand{\R}{I\!\! R}
\newcommand{\N}{I\!\! N}
\newcommand{\Z}{I\!\! Z}
\newcommand{\C}{I\!\! C}

\begin{document}

\bigskip
\bigskip

 Aristide Tsemo,

Department of Mathematics and Computer Sciences

Ryerson University,

350 Victoria Street, Toronto, ON

M5B 2K3

tsemo58@yahoo.ca

 Isaac Woungang,

Department of Mathematics and Computer Sciences

Ryerson University,

350 Victoria Street, Toronto, ON  M5B 2K3

iwoungan@scs.ryerson.ca
\bigskip
\bigskip

\centerline{\bf An index Theorem for gerbes.}

\bigskip

\centerline{\bf Abstract.}

{\it In this paper, we define the notion of vectorial gerbe, an
example of a vectorial gerbe is the Clifford-gerbe defined on a
riemannian manifold. We  show an index theorem for these objects.
This enables us to  study  the topology of riemannian manifolds.}

\bigskip

\centerline{\bf Introduction.}

\bigskip

In differential topology, the curvature tensors  of riemannian
metrics are very useful tools, for example, a  compact manifold
endowed with a metric whose  sectional curvature is constant is
the quotient of one of the following space by a discrete subgroup
of isometries: a vector space endowed with the flat riemannian
metric, the simply connected hyperbolic space, or the sphere
$S^n$.  If the  second Stiefel-Whitney class of a manifold $N$
vanishes, there exists a vector bundle called the bundle of
spinors which  has many applications in differential topology.
Suppose that the second Stiefel-Whitney class $w_2(N)$ of $N$ is
not zero, we deduce from the Giraud classification theorem, the
existence of a gerbe whose classifying cocycle is $w_2(N)$ which
is  according to Brylinski, and McLaughlin[2],  an illuminating
example of gerbe.   It is natural to study the relations between
this gerbe and the topology of $N$. For example one may expect to
generalize the Lichnerowicz theorem for spinors. On this purpose,
we need first to prove an index theorem for gerbes, which is our
purpose.

In their preprint [6], Murray and Singer have shown an index
theorem for bundle gerbes which is an example of the index formula
of this paper.

\bigskip

{\bf 1. On the notion of vectorial gerbes.}

\bigskip

The aim of this section is to develop the notion of vectorial
gerbes. First, we recall the definition of the  notion of sheaf of
categories on manifolds.

\medskip

{\bf  Definition 1.}

 Let $N$ be a manifold, a sheaf  of categories $C$ on $N$, is
a map $U\rightarrow C(U)$, where $U$ is an open subset of $N$, and
$C(U)$ a category which satisfies the following properties:

- For each inclusion $U\rightarrow V$, there exists a functor
$r_{U,V}:C(V)\rightarrow C(U)$ such that $r_{U,V}\circ
r_{V,W}=r_{U,W}$.

- Gluing conditions for objects,

Consider  an open covering family $(U_i)_{i\in I}$ of an open
subset $U$ of $N$, and for each $i$, an object $x_i$ of $C(U_i)$,
suppose that there exists a map $g_{ij}:r_{U_i\cap
U_j,U_j}(x_j)\rightarrow r_{U_i\cap U_j, U_i}(x_i)$. Denote by
$g_{i_1i_2}^{i_3}$, the restriction of the map $g_{i_1i_2}$ to
$U_{i_1}\cap U_{i_2}\cap U_{i_3}$, and suppose that
${g_{i_1i_2}}^{i_3}{g_{i_2i_3}}^{i_1}={g_{i_1i_3}}^{i_2}$, then
there exists an object $x$ of $C(U)$ such that $r_{U_i,U}(x)=x_i.$

- Gluing conditions for arrows,

Consider two objects $P$ and $Q$ of $C(N)$,  the map $U\rightarrow
Hom(r_{U,N}(P),r_{U.N}(Q))$ defined on the category of open
subsets of $N$ is a sheaf.

Moreover, if the following conditions are satisfied, the sheaf of
categories $C$ is called a gerbe.

$G1$

There exists an open covering family $(U_i)_{i\in I}$ of $N$, such
that for each $i$, the category $C(U_i)$ is not empty.

$G2$

Let $U$ be an open subset of $N$, for each objects $x$ and $y$ of
$C(U)$, there exists an open covering family $(U_i)_{i\in I}$ of
$U$, such that $r_{U_i,U}(x)$ and $r_{U_i,U}(y)$ are isomorphic.

$G3$

Every arrow of $C(U)$ is invertible. There exists a sheaf $L$ of
sections of a principal bundle over $N$ such that for each object
$x$ of $C(U)$, $Hom(x,x)$ is isomorphic to $L(U)$, by an
isomorphism which commutes with restriction maps.

The sheaf $L$ is called the band of the gerbe $C$, in the sequel,
we will consider only gerbes with commutative band.

\bigskip

{\bf Notations.}

\bigskip

 For a covering family $(U_i)_{i\in I}$ of $N$, and an object $x_{i_1}$
 of $C(U_{i_1})$, we denote
by ${x_{i_1}}^{i_{i_2..i_n}}$ the element $r_{U_{i_1}\cap..\cap
U_{i_n},U_{i_1}}(x_{i_1})$,  by $U_{i_1..i_n}$ the intersection
$U_{i_1}\cap..\cap U_{i_n}$, and for a map
$h_{i_1}:e_{i_1}\rightarrow e'_{i_1}$ between two objects,
$e_{i_1}$ and $e'_{i_1}$ of $C(U_{i_1})$, we denote by
${h_{i_1}}^{i_2..i_p}$, its restriction to $U_{i_1..i_p}$. Suppose
that the objects of the category $C(U)$ are vectors bundles, and
consider  $s_{i_1}:U_{i_1}\rightarrow e_{U_{i_1}}$,  a section of
the bundle $e_{U_{i_1}}$ defined over $U_{i_1}$, we denote by
$s_{i_1}^{i_2...i_n}$ its restriction to $U_{i_1..i_n}$.

\bigskip

{\bf Definition 2.}

A gerbe is a vectorial gerbe, if and only if for each open subset
$U$ of $N$, the category $C(U)$ is a category whose objects are
vector bundles over $U$, which typical fiber is the vector space
$V$, and such that the maps between objects of $C(U)$ are
isomorphisms of vector bundles. The vector space $V$ is called the
typical fiber of the vectorial gerbe. The classifying cocycle of
the gerbe is defined as follows: there exists an open covering
family $(U_i)_{i\in I}$ of $N$, a commutative subgroup $H$ of
$Gl(V)$, such that there exist maps ${g'}_{ij}:U_i\cap
U_j\longrightarrow Gl(V)$, which define isomorphisms
$$
g_{ij}: U_i\cap U_j\times V\longrightarrow U_i\cap U_j\times V
$$
$$
(x,y)\longrightarrow (x,g_{ij}(x)y)
$$
and such that
$c_{i_1i_2i_3}={g_{i_1i_2}}^{i_3}{g_{i_2i_3}}^{i_1}{g_{i_3i_1}}^{i_2}$
is the classifying Cech $2$-cocycle which takes its values in the
sheaf of $H$-valued functions.

\medskip

{\bf Example.}

\medskip

{\bf The Clifford gerbe associated to a riemannian structure.}

\medskip

Let $N$ be a riemannian $n$-manifold, and $O(N)$ the reduction of
the bundle of linear frames which defines the riemannian structure
of the manifold $N$. The bundle $O(N)$ is a locally trivial
principal bundle over $N$ which typical fiber is $O(n)$, the
orthogonal group of $n\times n$ matrices. We have the exact
sequence $1\rightarrow {{\Z}/2}\rightarrow Spin(n)\rightarrow
O(n)\rightarrow 1$, where $Spin(n)$ is the universal cover of
$O(n)$. We can associate to this problem a gerbe whose band is
${{\Z}/2}$, defined as follows:
 for each open subset $U$ of $N$, $Spin(U)$ is the category of
$Spin(n)$-bundles over $U$ such that the quotient of each object
of $Spin(U)$  by ${{\Z}/2}$  is the restriction $O(U)$, of $O(N)$
to $U$. The classifying $2$-cocycle of this gerbe is the second
Stiefel-Whitney class.

One can associate to this gerbe, a vectorial gerbe named the
Clifford gerbe $Cl(N)$ which is defined as follows: for each open
subset $U$ of $N$, $Cl(U)$ is the category whose objects are
Clifford bundles  associated to the objects of $Spin(U)$. The
gerbe $Cl(N)$ is a vectorial gerbe.

Let ${g'}_{ij}$ be the transitions functions of the bundle $O(N)$,
for each map $g_{ij}$, consider an element ${g}_{ij}$ over
${g'}_{ij}$ which take its values in $Spin(n)$. Then the element
$g_{ij}(x)$ acts on $Cl({\R}^n)$ by left multiplication, we will
denote by $h_{ij}(x)$ the resulting automorphism of $Cl({\R}^n)$.
The classifying cocycle of the Clifford gerbe is
 defined by the maps ${h_{i_1i_2}}^{i_3}{h_{i_2i_3}}^{i_1}{h_{i_3i_1}}^{i_2}$.

\bigskip

{\bf The gerbe defined by  the lifting problem associated to a
vectorial bundle.}

\medskip

Consider a vector bundle $E$ over $N$, whose typical fiber is the
vector space $V$. One classically associate to $E$, a principal
$Gl(V)$-bundle. We suppose that this bundle has a reduction $E_L$,
where $L$ is a subgroup of $Gl(V)$. Consider a central extension
$1\rightarrow H\rightarrow G\rightarrow L\rightarrow 1$. We denote
by $\pi:G\rightarrow L$, the projection and suppose it has local
sections. This central extension defines a gerbe $C_H$ on $N$,
such that for each open subset $U$ of $N$, the objects of $C_H(U)$
are $G$-principal bundles over $U$, whose quotient by $H$ is the
restriction of $E_L$ to $U$.

Suppose moreover defined a representation $r:G\rightarrow Gl(W)$,
and a surjection $f:W\rightarrow V$ such that the following square
is commutative:
$$
\matrix{W &{\buildrel{r(h)}\over{\longrightarrow}}& W\cr
\downarrow f &\\& f \downarrow\cr V&{\buildrel{\pi(h)}\over
{\longrightarrow}} &V}\leqno (1)
$$

then one can define the vectorial gerbe $C_{H,W}$ on $N$ such that
the objects of $C_{H,W}(U)$ are fiber products $e_U\propto r$,
where $e_U$ is an object of $C_H(U)$.  Let $(U_i)_{i\in I}$ be a
trivialization of $E$ defined by the transitions functions
${g'}_{ij}$, we consider a map $g_{ij}:U_i\cap U_j\rightarrow G$
over ${g'}_{ij}$. The classifying cocycle of the gerbe $C_{H,W}$
is defined by
$r({g_{i_1i_2}}^{i_3})r({g_{i_2i_3}}^{i_1})r({g_{i_3i_1}}^{i_2})$.

\medskip

{\bf Definition 3.}

-A riemannian metric on a vectorial gerbe $C$ is defined by the
following data:

For each object $e_U$ of $C(U)$, a riemannian metric $<,>_{e_U}$
on the vector bundle $e_U$. We suppose that  each morphism
$h:e_U\rightarrow e'_U$, between the objects $e_U$ and $e'_U$ of
$C(U)$, is an isomorphism between the riemannian manifolds
$(e_U,<,>_{e_U})$ and $(e'_U,<,>_{e'_U})$. We remark that the band
need to be contained in a sheaf of sections of a principal bundle
whose typical fiber is a compact group in this case, since its
elements preserve the riemannian metric.

\medskip

An example of a scalar product on a gerbe is the following:
consider the Clifford gerbe $Cl(N)$, we know that the group
$Spin(n)$ is a compact group, its action on $Cl({\R}^n)$ preserves
a scalar product. This scalar product defines on each fiber of an
object $e_U$ of $Cl(U)$, a scalar product which defines a
riemannian metric $<,>_{e_U}$ on $e_U$. The family of riemannian
metrics $<,>_{e_U}$ is a riemannian metric defined on the gerbe
$Cl(N)$.

\medskip

{\bf Definition 4.}

 - A global section of a vectorial gerbe $C$ associated to a $1$-Cech
chain $(g_{ij})_{ij}$, is defined by an open covering family
$(U_i)_{i\in I}$ of $N$, for each element $i$ of $I$, an object
$e_i$ of $C(U_i)$, a section $s_i$ of $e_i$, a family of morphisms
$g_{ij}:e^i_j\rightarrow e^j_i$ such that on $U_{ij}$ we have
$s_i=g_{ij}(s_j)$.

Let $(s_i)_{i\in I}$ be a global section, on $U_{i_1i_2i_3}$ we
have: ${s_{i_2}}^{i_3}=g_{i_2i_3}({s_{i_3}}^{i_2})$,
${s_{i_1}}^{i_2}=g_{i_1i_2}({s_{i_2}}^{i_1})$. This implies that
the restriction ${s_{i_1}^{i_2i_3}}$ of $s_{i_1}$ to
$U_{i_1i_2i_3}$ verifies:
${s_{i_1}}^{i_2i_3}={g_{i_1i_2}}^{i_3}{g_{i_2i_3}}^{i_1}({s_{i_3}}^{i_1i_2})={g_{i_1i_3}}^{i_2}({s_{i_3}}^{i_1i_2})$.
Or equivalently
${{g_{i_1i_3}}^{i_2}}^{-1}{g_{i_1i_2}}^{i_3}{g_{i_2i_3}}^{i_1}(s_{i_3})=s_{i_3}$.
Remark that the restriction of the element $s_{i_3}$ to
$U_{i_1i_2i_3}$ is not necessarily preserved by all the band.

Suppose that $N$ is compact, and $I$ is finite. We can suppose
that there exists $i_0$ such that $T=U_{i_0}-\cup_{i\neq i_0}U_i$
is not empty. Consider a section $s_{i_0}$ of $e_{i_0}$ whose
support is contained in $T$, then we can define a global section
$(u_i)_{i\in I}$ such that $u_{i_0}=s_{i_0}$, and if $i\neq i_0$,
$u_i=0$. This ensures that $S(g_{ij})$ is not empty. We will
denote by $S(g_{ij})$ the family of global sections associated to
$(g_{ij})_{i,j\in I}$. Remark that $S(g_{ij})$ is a vector space.

\bigskip
\medskip

{\bf Proposition 5.}

{\it Suppose that the vectorial gerbe $C$ is the gerbe associated
to the lifting problem defined by the extension $1\rightarrow
H\rightarrow G\rightarrow L\rightarrow 1$ and the vector bundle
$E$. Let $r:G\rightarrow Gl(W)$ be a representation, we suppose
that the condition  of the  diagram (1) is satisfied. Then for
each  $G$-chain $g_{ij}$, each element $(s_i)_{i\in I}$ of the
vector space of global sections $S(g_{ij})$, there exists a
section $s$ of $E$, such that $s_{\mid U_i}=f\circ s_i$.}

\medskip

{\bf Proof.}

Let $(s_i)_{i\in I}$ be a global section associated to the chain
$S(g_{ij})$, then on $U_{ij}$, we have $s^j_i=g_{ij}(s^i_j)$,
(where $s^j_i$ is the restriction of $s_i$ to $U_{ij}$) this
implies that on $U_{ij}$, $f(s^j_i)=g'_{ij}(f(s^i_j))$. Thus the
family $(f(s_i))_{i\in I}$ of local sections of $E$ defines a
global section $s$ of $E$.

\medskip

{\bf Remark.}

\medskip

Let $s$ be a section of the bundle $E$, locally we can define a
family of sections $s_i$ of $e_i$, such that $f(s_i)=s_{\mid
U_i}$. We can consider the chain $s_{ij}=s^j_i-g_{ij}(s^i_j)$. The
family ${s_{jk}}-s_{ik}+s_{ij}$ is a $2$-cocycle. Whenever there
exists a chain $g_{ij}$, a global section $s=(s_i)_{i\in I}$ such
that $s_i=g_{ij}(s_j)$, and $f(s_i)=s_{\mid U_i}$, it is not sure
that such a global section exists for another chain $h_{ij}$. This
motivates the following definition:

\medskip

{\bf Definition 6.}

We  define the vector space $L$ of formal global sections of the
vector gerbe $C$, to be the vector space  generated by the
elements $[s]$, where $s$ is an element of a set of global
sections $S(g_{ij})$. The elements of $L$, are  finite sum of
global sections.

\bigskip

{\bf The Prehilbertian structure of $S(g_{ij})$.}

\medskip

First we remark that $S(g_{ij})$ is a vector space. Let $s$, and
$t$ be elements of $S(g_{ij})$, we will denote by $s_i$ and $t_i$
the sections of $e_i$ which define respectively the global
sections $s$ and $t$. We have $s^j_i=g_{ij}(s^i_j)$ and
$t^j_i=g_{ij}(t^i_j)$ this implies that
$as^j_i+bt^j_i=g_{ij}(as^i_j+bt^i_j)$, where $a$ and $b$ are real
numbers.

\medskip

- The scalar structure of $S(g_{ij})$.

\medskip

Let $(V_k, f_k)_{k\in K}$ be a partition of unity subordinate to
$(U_i)_{i\in I}$, this implies that for each $k$ there exists an
$i(k)$ such that $V_k$ is a subset of $U_{i(k)}$. Since the
support of $f_k$ is a compact subset of $V_k$, we can calculate
$\int_{V_k}f_k<s_{ii(k)},t_{ii(k)}>$ where $s_{ii(k)}$ and
$t_{ii(k)}$ are the respective restrictions of $s_{i(k)}$ and
$t_{i(k)}$ to $V_k$. Remark that since we have supposed that
$s^j_i=g_{ij}(s^i_j)$ and $g_{ij}$ is a riemannian isomorphism
between $e_j^i$ and $e_i^j$, if $V_k$ is also included in
$U_{j(k)}$, then we have the equality $<s_{ii(k)},t_{ii(k)}>=
<s_{jj(k)},t_{jj(k)}>$ on $V_k$. We can define
$<s,t>=\sum_k\int<f_ks_{ii(k)},f_kt_{ii(k)}>$.

We will denote by $L ^2(S(g_{ij}))$ the Hilbert completion of the
Pre Hilbert structure of $(S(g_{ij}),<,>)$.

\bigskip

{\bf The scalar structure on the set of formal global sections
$L$.}

\medskip

Let $s$ and $t$ be two formal global sections, we have
$s=[s_{n_1}]+..+[s_{n_p}]$, and $t=[t_{m_1}]+..+[t_{m_q}]$, where
$s_{n_i}$ and $t_{m_j}$ are global sections.

We will define a scalar product on $L$ as follows: if $s$ and $t$
are elements of the same set of global sections $S(g_{ij})$,
$<[s],[t]>=<s,t>_{S(g_{ij})}$, where $<s,t>_{S(g_{ij})}$ is the
scalar product defined at the paragraph above. If $s$ and $t$ are
not elements of the same set of global sections, then
$<[s],[t]>=0$.

\medskip

{\bf Proposition 7.}

{\it An element of $L^2(S(g_{ij}))$ is a family of $L^2$ sections
$s_i$ of $e_i$ such that $s^j_i=g_{ij}(s^i_j)$.}

\medskip

{\bf Proof.}

Let $(s^l)_{l\in N}$ be a Cauchy sequence of $(S(g_{ij}),<,>)$. We
can suppose that the open sets $V_k$ used to construct the
riemannian metric are such that the restriction $e^k_i$, of $e_i$,
to $V_k$ is a trivial vector bundle. The sequence
$(f_ks^l_{ii(k)})$ is a Cauchy sequence defined  on the support
$T_k$ of $f_k$. Since this support is compact, we obtain that
$(f_ks^l_{ii(k)})_{l\in N}$ goes to an $L^2$ section $s_{i(k)}$ of
$e^k_i$. We can define $s_i=\sum_{k,V_k\cap
U_i\neq\phi}f_ks_{i(k)}$. The family $(s_i)_{i\in I}$ defines the
requested limit.

\medskip

Suppose that morphisms between objects commute with the laplacian
$\Delta$, we can then endow $S(g_{ij})$ with the prehilbertian
structure defined by $<u,v>=\int<\Delta^s(u),v>$, where $s$ a
positive real number, and $\Delta^s(u)$ is the global section
defined by ${\Delta^s(u)}_i=\Delta^s(u_i)$. We will denote by
$H_s(S(g_{ij}))$ the Hilbert completion of this prehilbertian
space.

We will define the formal $s$-distributional global sections
$H_s(L)$ as the vector space generated by finite sums
$[s_1]+..[s_k]$ where $s_i$ is an element of an Hilbert space
$H_s(S(g_{ij}))$.

\bigskip

{\bf Connection on riemannian gerbes and characteristic classes.}

\medskip

The notion of connection is not well-defined for general vectorial
gerbes, nevertheless the existence of a riemannian structure on a
vectorial gerbe $C$, defined on the manifold $N$, gives rise to a
riemannian connection on each object $e_U$ of $C(U)$, this family
of riemannian connections will be the riemannian connection of the
gerbe $C$.

Let $(U_p)_{p\in P}$ be an open covering of $N$, and $so(V)$ the
Lie-algebra of the orthogonal group $SO(V)$ of $V$. Suppose that
the objects of $C(U_p)$ are trivial bundles. The riemannian
connection of the object $e_p$ of $C(U_p)$ is defined by a $so(V)$
$1$-form $w_p$ on $TU_p$, and the covariant derivative of this
connection evaluated at a section $s_p$ of $e_p$ is $ds_p+w_ps_p$.
The curvature of this connection is the $2$-form
$\Omega_p=dw_p+w_p\wedge w_p$.

The $2k^{eme}$ Chern form of $e_p$, ${c^p}_{2k}$ is defined by the
trace $Trace[({i\over 2\pi}\Omega_p)^k]$. Let $e'_p$ be another
object of $C(U_p)$. There exist isomorphisms
$\phi_p:e_p\rightarrow U_p\times V$,  $\phi_p':e'_p\rightarrow
U_p\times V$, and $g_p:e_p\rightarrow e'_p$. The map $\phi_p'\circ
g_p\circ {\phi_p}^{-1}$ is an automorphism of $U_p\times V$
defined by the map:
 $h_p:(x,y)\rightarrow (x,u_p(x)y))$.
The riemannian ${\phi_p}^{-1}<,>={\phi'_p}^{-1}<,>$ is preserved
by $h_p$. This implies that ${c^p}_{2k}$ is equal to the  Chern
$2k$-form of $e'_p$.

\bigskip

 There exists an isomorphism between the respective restrictions
 $e^p_j$ of $e_j$, and $e^j_p$ of $e_p$ to $U_{pj}$. As above we can show that this implies
 that the Chern $2k$-forms of $e^p_j$ and $e^j_p$ coincide  on $U_p\cap
 U_j$. We can define a global form $c_{2k}(N)$ on $N$ whose restriction to $U_p$ is $c^p_{2k}$,
  which will be called the $2k$-Chern
 form of the riemannian gerbe.

 We can define the $c(C)=c_1(N)+...+c_n(N)$ the total Chern form
 of the gerbe, and the form defined locally by ${ch(C)}_{\mid U_p}=
 Trace(exp(i{\Omega\over 2\pi}))$ the total Chern character.

\medskip

{\bf 2. Operators on riemannian gerbe.}

\medskip

We begin by recalling the definition of pseudo-differential
operators for open subsets of ${\R}^n$, and for manifolds.
 Let $U$, be an open subset of
${\R}^n$, we denote by $S^m(U)$ the set of smooth functions
$p(x,u)$ defined on $U\times{\R}^n$, such that for every compact
set $K\subset U$, and every multi-indices $\alpha$ and $\beta$, we
have $\mid\mid
D^{\alpha}D^{\beta}p(x,u)\mid\mid<C_{\alpha,\beta,K}(1+\mid\mid
u\mid\mid)^{m-\mid \alpha\mid}$, where $D^{\alpha}$ is the partial
derivative in respect to $\alpha$.

Let $K(U)$ and $L(U)$ denote respectively the space of smooth
functions with compact support defined on $U$, and the space of
smooth functions on $U$. We can define the map $P:K(U)\rightarrow
L(U)$ by:
$$
P(f)=\int_U p(x,u)\hat f(u)e^{i<x,u>}du
$$
where $\hat f$ is the Fourier transform of $f$.

\medskip

{\bf Definition 1.}

 An operator on $U$  is pseudo-differential, if it is locally of
the above type.

\medskip

{\bf Definition 2.}

Let $P$ be a pseudo-differential operator, $(U_i)_{i\in I}$ an
open covering family of $U$ such that the restriction of $P$ to
$U_i$ is defined by $P(f)=\int_{U_i} p_i(x,u)\hat
f(u)e^{i<x,u>}du$. The operator is of degree $m$ if
${\sigma(p)}_{\mid{U_i}}=lim_{t\rightarrow\infty}{p_{i}(x,tu)\over
t^m}$ exists. In this case $\sigma(p)$ whose restriction to $U_i$
is ${\sigma(p)}_{U_i}$ is called the symbol of $P$.

\medskip

Let $E$ be a vector bundle over the riemannian manifold $N$,
endowed with a scalar metric. We denote by $K(E)$ and $L(E)$ the
respective  space of smooth sections of $E$ with compact support,
and the space of smooth sections of $E$. An operator on the vector
bundle $E$, is a map $P:K(E)\rightarrow L(E)$ such that there
exists an open covering family $(U_i)_{i\in I}$ of $N$ which
satisfies:

- The restriction of $E$ to $U_i$ is trivial, in fact we suppose
that $U_i$ is the domain of a chart.

- The restriction of $P_i$, of $P$ to $U_i$ is a map
$P_i:K(U_i\times V)\rightarrow L(U'_i\times V)$ where $V$ is the
typical fiber of $E$, and $U'_i$ an open chart of $N$.

- If we consider charts $\phi_i$ and $\psi_i$ whose domains are
respectively $U_i$ and $U'_i$, and such that $\phi_i(U_i\times
V)=\psi_i(U'_i\times V)=U\times {\R}^n$, then the map $P_i$ is
defined by  a matrix $(p_{kl})$, where $p_{kl}$
 defines an operator of degree $m$. More precisely, if $s'$ is a section
  of $E$ over $U_i$,
and $s=(s_1,..,s_n)=\phi_i(s')$, we can define
$t_k=\sum_{l=1}^{l=n}\int p_{kl}(x,u)\hat s_l(u)e^{i<x,u>}du$, and
$P_i(s')={\psi_i}^{-1}(t_1,..,t_n)$.

 Consider $SN$ the sphere bundle of the cotangent space $\pi:T^*N\rightarrow N$,
and $\pi^*E$ the pull-back of $E$ to $T^*N$, the symbols defined
by the operators $p_{ij}$ define a map $\sigma:\pi^*E\rightarrow
\pi^*E$. Consider now the projection $\pi_S:SN\rightarrow N$, then
$\sigma$ induces a map $\sigma_S:{\pi_S}^*E\rightarrow
{\pi_S}^*E$.

 Let $s$ be a positive integer
 we denote by ${{H_s}^{loc}}(N,E)$, the space of distributional sections $u$ of $E$
such that $D(u)$ is a ${L^2}$-section, where $D$ is any
differential operator of order less than $s$, and by
${{H_s}^{comp}}(N,E)$ the subset of elements of
${{H_s}^{loc}}(N,E)$ with compact support. Remark that if $N$ is
compact, then ${{H_s}^{loc}}(N,E)={{H_s}^{comp}}(N,E)$ in this
case, we denote this space by $H_s(N,E)$. We define by
${H_{-s}}^{loc}(N,E)$ to be the dual space of ${H_s}^{comp}(N,E)$,
and by ${H_{-s}}^{comp}(N,E)$ the dual space of
${H_s}^{loc}(N,E)$.

Suppose that $N$ is compact, then the Sobolev space $H_s(N,E)$ is
a Hilbert space endowed with the norm defined by $(\mid\mid \int_N
<\Delta^su,u>\mid\mid)^{1\over 2}$. Every operator $P$ of order
less than $m$ can be extended to a continuous morphism
$H_s(N,E)\rightarrow H_{s-m}(N,E)$.

\medskip

We will adapt now the definition of an operator on a manifold to
the definition of operators on vectorial gerbes.

\medskip

{\bf Definition 3.}

Let $C$ be a riemannian gerbe defined on the manifold $N$, an
operator $D$ of degree $m$ on $C$, is a family of operators
$D_{e}$ of degree $m$ defined on each oject $e$,   of the category
$C(U)$, for each open subset $U$ of $N$. We suppose that for each
morphism $g:e\rightarrow f$, $D_fg^*=g^*D_e$, where $g^*(s)=g(s)$.

\medskip

{\bf Remark.}

The last condition in the previous definition  implies that $D_e$
is invariant by the automorphisms of $e$. The operators considered
in the sequel will be assumed to be continue, and we will assume
that they preserve $C^{\infty}$-sections. Thus an operator defines
a map $D^H_e:{{H_s}^{comp}}(U,e)\rightarrow {H_{s-m}^{loc}}(U,e)$.

\medskip

{\bf Proposition 4.} {\it Let $D$ be an operator of degree $n$
defined on the riemannian gerbe $C$, then $D$ induces a map
$D_{S(g_{ij})}:H_s(S(g_{ij}))\rightarrow H_{s-n}(S(g_{ij}))$ and a
map $D_L:H_s(L)\rightarrow H_{s-n}(L)$.}

\medskip

{\bf Proof.}

Consider a global distributional section $s$ which is an element
of $H_s(S(g_{ij}))$, we have
$g_{ij}(D_{e^i_j}(s^i_j))=D_{e^j_i}(s^j_i)$. This implies the
result.

\bigskip

{\bf The symbol of an operator defined on a gerbe.}

\bigskip

Let $C$ be a vectorial gerbe defined on a compact manifold $N$
endowed with the operator $D$ of degree $m$, for each object $e$
of $C(U)$, we can pull back the bundle $e$ by the projection map
$\pi_{SU}:SU\rightarrow U$ to a bundle $\pi_{SU}^*e$ over $SU$,
where $SU$ is the restriction of the cosphere bundle defined by a
fixed riemannian metric of $T^*N$. The set $C_S(U)$ which elements
are $\pi_{SU}^*e$ where $e$ is an object of $C(U)$ is a category.
The maps between objects of this category are induced by maps
between elements of $C(U)$. The map $U\rightarrow C_S(U)$ is a
gerbe which has the same band than $C$. Now on the object $e$, we
can define a symbol $\sigma_{D_e}:\pi_{SU}^*e\rightarrow
\pi_{SU}^*e$.

Remark that for every automorphism $g$ of $e$, the fact that
$g^*\circ D_e=D_e\circ g^*$ implies that
$\sigma_{gD_eg^{-1}}=\sigma_{D_e}$.

\bigskip
\medskip

{\bf Proposition 5. Rellich Lemma for the family $S(g_{ij})$.}

{\it Let $(f_n)_{n\in {\N}}$ be a sequence of elements of
$H_s(S(g_{ij}))$, we suppose that there exists a constant $L$ such
that $\mid\mid f_n\mid\mid_s<L$, then for every $s>t$,  there
exists a subsequence $f_{n_k}$ which converges in $H_t$.}

\medskip

{\bf Proof.}

Let $(s_n)_{n\in {\N}}$ be a sequence of sections which satisfy
the condition of the proposition, and $(V_{\alpha}, f_{\alpha})$ a
partition of unity subordinate to $(U_i)_{i\in I}$. We suppose
that the support of $f_{\alpha}$ is a compact space $K_{\alpha}$.
We denote by $s^i_n$ the section of $e_i$ which defines $s_n$, and
by $s^i_{n\alpha}$ the restriction of $s^i_n$ to the restriction
of $e_i$ to $V_{\alpha}$. The family $(f_{\alpha}s^i_{n\alpha})$
$H_t$-converges  towards the element $s_{i\alpha}$ in $V_{\alpha}$
by the classical Rellich lemma. We can write then
$s^i=\sum_{V_{\alpha}\cap U_i}f_{\alpha}s_{i\alpha}$. this is an
$H_t$ map since the family of $V_{\alpha}$ can be supposed to be
finite, since $N$ is compact. The family $(s^i)$ defines a global
$H_t$-section which is the requested limit.

\medskip

{\bf Remark.}

A compact operator between Hilbert spaces, is an operator which
transforms bounded spaces to compact spaces. The previous lemma
implies that if $s>t$, then the inclusion
$H_s(S(g_{ij}))\rightarrow H_t(S(g_{ij}))$ is compact, indeed
since $H_s(S(g_{ij}))$ is a separate space, a compact subspace of
$H_s(S(g_{ij}))$ is a set such that we can extract a convergent
sequence from every  sequence.

\medskip

{\bf Proposition 6.}

{\it The space $Op(C)$  of continuous linear maps of
$H_s(S(g_{ij}))$ is a Banach space.}

\medskip

{\bf Proof.}

Let $(D_n)_{n\in {\N}}$ be a Cauchy sequence of elements of
$Op(C)$, for each global section $s$, the sequence $D_n(s)$  is a
Cauchy sequence in respect to the norm of $H_s(S(g_{ij}))$, we
conclude that it converges towards an element $D(s)$. The map
$D:s\rightarrow D(s)$ is the requested limit. It is bounded since
$(D_n)_{n\in N}$ is a Cauchy sequence.

\medskip

The previous proposition allows us to define $O^n$, the completion
of the pseudo-differential operators in   $OP(C)$ of order $n$,
and to extend the symbol $\sigma$ to $O^n$. Now we will show that
the kernel of the extension of the symbol to $O^n$ contains only
compact operators.

\medskip

{\bf Proposition 7.}

{\it The kernel of the symbol map contains only compact
operators.}

\medskip

{\bf Proof.}

The symbol $\sigma(P)$ of the operator $P$ is zero if and only if
the order $m$ of the operator is  less than $-1$. This implies
that the operator $P$ is compact. To see this, we can suppose our
operator to be an $L^2(S(g_{ij}))$ operator, by composing it by
the inclusion map $H_{2-m}(S(g_{ij}))\rightarrow L^2(S(g_{ij})$,
we conclude by using the previous Rellich lemma.

\medskip

We end this section by defining elliptic operators for gerbes.

\medskip

{\bf Definition 8.}

We say that an operator is elliptic if the family of symbols
$\sigma_{D_e}$ are invertible maps.

\medskip

{\bf 3. $K$-theory and the index.}

\medskip

In this part we will give the definitions of the $K$-theory groups
$K_0$, and $K_1$, and show how we can use them to associate to a
symbol of an operator on a gerbe defined on $N$, an element of
$K_0(T^*N)$.

\medskip

We will denote by $M_n$ the vector space of $n\times n$ complex
matrices. For $n\leq m$, we consider the natural injection
$M_n\rightarrow M_m$. We will call $M_{\infty}$, the inductive
limit of the vector spaces $M_n, n\in{\N}$.

Let $R$ be a ring,  $p$, and $q$ be two idempotents of
$R_{\infty}=R\otimes M_{\infty}$, we will say that $p\sim q$ if
and only if there exists elements $u$ and $v$ of $R_{\infty}$ such
that $p=uv$, and $q=vu$. We denote by $[p]$ the class of $p$, and
by $Idem(R_{\infty})$ the set of equivalence classes.
 If $[p]$ and $[q]$ are represented respectively by elements of $R\otimes M_n$ and
$R\otimes M_m$, we can define an idempotent of $R\otimes M_{n+m}$
represented by the matrix $\pmatrix{p&0\cr 0&q}$ that we denote
$[p+q].$

\medskip

{\bf Definition 1.}

We will denote by $K_0(R)$, the semi-group $Idem(R_{\infty})$,
endowed with the law $[p]+[q]=[p+q]$.

\medskip

Let $N$ be a compact manifold, and $C(N)$ the set of complex
valued functions on $N$. It is a well-known fact that for a
complex vector bundle $V$ on $N$, there exists  a bundle $W$, such
that $V\oplus W$ is a trivial bundle  isomorphic to $N\times
{\C}^l$. We can thus identify a vector bundle over $N$, to an
idempotent of $C(X)\otimes M_l$ which is also an idempotent of
$C(X)_{\infty}$. This enables to identify $K_0(N)$ to $K_0(C(X))$.
In fact the semi-group $K_0(C(X))$ is a group.

\medskip

Let $Gl_n(R)$ be the group of invertible elements of $M_n(R)$, if
$l\leq n$ we have the canonical inclusion map $Gl_l(R)\rightarrow
Gl_n(R)$. We will denote by $Gl_{\infty}(R)$, the inductive limit
of the groups $Gl_n(R)$.

\medskip

{\bf Definition 2.}

Let $Gl_{\infty}(R)_0$ be the connected component of
$Gl_{\infty}(R)$. We will denote by $K_1(R)$ the quotient of
$Gl_{\infty}(R)$ by $Gl_{\infty}(R)_0$.

\medskip

For a compact manifold $N$, we denote  $K_1(C(N))$ by $K_1(N)$.

Consider now an exact sequence $0\rightarrow R_1\rightarrow
R_2\rightarrow R_3\rightarrow 0$ of $C^*$-algebras, we have the
following exact sequence in $K$-theory:
$$
K_1(R_1)\rightarrow K_1(R_2)\rightarrow K_1(R_3)\rightarrow
K_0(R_1)\rightarrow K_0(R_2)\rightarrow K_0(R_3).
$$

\medskip

Let ${\cal H}$ be an Hilbert space, we denote by $B({\cal H})$ the
space of continuous operators defined on ${\cal H}$, and ${\cal
K}$ the subspace of compact continuous operators of $B({\cal H})$.
We have the following exact sequence:
$$
0\rightarrow {\cal K}\rightarrow B({\cal H})\rightarrow B({\cal
H})/{\cal K}=Ca\rightarrow 0.
$$

The algebra $Ca$ is called the Calkin algebra of ${\cal H}$. It is
a well-known fact that $K_0({\cal K})={\Z}$.

\medskip

 Let $N$ be a riemannian manifold, and $C$ a riemannian gerbe
defined on $N$. Consider an elliptic operator $D$ defined on $C$,
of degree $l$. The operator $D$ induces a morphism:
$D:L^2(S(g_{ij}))\rightarrow H_{2-l}(S(g_{ij}))$. Consider the
operator  $(1-\Delta)^{-m}$ of degree $-l$. The operator
$D'=(1-\Delta)^{-m}D$ is a morphism of $L^2(S(g_{ij}))$.  The
symbol of $(1-\Delta)^{-m}D$ is also $\sigma(D)$. This implies
that the image of the operator $D$ in the Calkin algebra of
$L^2(S(g_{ij}))$ is invertible. It thus defines a class
$[\sigma(D')]$ of $K_1(Ca)$. The image of $[\sigma(D')]$ in
 $K_0({\cal K})$ is the index of $D$. We remark that
the index of the operator depends only of the symbol.

For every object $e$ of $C(U)$, the symbol $\sigma({D_e})$ is an
automorphism of ${\pi_{SU}}^*e$, it defines an element
$[\sigma(D_e)]$ of $K_1(C(S(U)))$ (recall that $S(U)$ is the
cosphere bundle over $U$ defined by the riemannian metric).

\bigskip

{\bf Remark.}

Let $U$ be an open subset such that the objects of $C(U)$ are
trivial bundles. Consider an object $e$ of $C(U)$, and a
trivialization map $\phi_e:e\rightarrow U\times V$. For every
object $f$ of $C(U)$, we have
${{\phi_e}^{-1}}^*(\sigma_{D_e})={{\phi_f}^{-1}}^*(\sigma_{D_f})$.

\medskip

{\bf Proposition 3.}

{\it Let $C$ be a vectorial gerbe defined over the compact
manifold $N$, and $U$ an open subset of $N$, there exists a
trivial complex bundle $f_n=N\times {\C}^n$, such that each object
of $C(U)$  is isomorphic to  a sub-bundle of the restriction of
$f_n$ to $U$.}

\medskip

{\bf Proof.}

Let $(U_i)_{i\in I}$ be an open  finite covering family of $N$
such that for each $i$ the objects of the category $C(U_i)$ are
trivial bundles. Let $e_U$ be an object of $C(U)$, the restriction
$e_i$ of $e$ to $U_i$ is a trivial vectorial bundle. We consider a
fixed object $e_i^0$ of $C(U_i)$. Consider a finite partition of
unity $(f_p)_{p\in 1,..,l}$ subordinate to the covering family
$(U_i)_{i\in I}$, and $g_i:e_i^0\rightarrow V$ the composition of
the trivialization and the second projection. Let
$h_i:e_i\rightarrow e_i^0$ be an isomorphism, we can define the
map $k:e_U\rightarrow {\C}^{ldimV} $ such that
$k(x)=(f_1(\pi_{e_U}(x))g_1h_1(x),..,f_l(\pi_{e_U}(x))g_lh_l(x))$.
$k$ induces a map $K:N\rightarrow G_p({\C}^n)$, where $p=dim V$
$n=ldim(V)$, and $G_p({\C}^n)$ is the Grassmanian of $p$-complex
plane of ${\C}^n$; $e_U$ is the pulls-back of the canonical
$p$-vector bundle over $G_{p}({\C}^n)$. We remark that it is a
subbundle of the pulls-back of the trivial bundle
$G_p({\C}^n)\times {\C}^n$.

\medskip

The orthogonal bundle of $e_U$ in the previous result can be
chosen canonically by considering the orthogonal bundle of the
canonical ${\C}^p$-bundle over $G_p({\C}^n)$ in $G_p({\C}^n)\times
{\C}^n$. We will suppose that this bundle is chosen canonically in
the sequel.

\medskip

 {\bf Proposition 4.}

{\it Let $C$ be a riemannian gerbe defined over the compact
manifold $N$. Then we can associate naturally to the symbol of the
elliptic operator $D$, a class $[\sigma_D]$ in $K_1(T^*N)$.}

\medskip

{\bf Proof.}

Let $(U_i)_{i\in I}$, be a finite open covering family of $N$ such
that for each $i$, each object of $e_i$ of $C(U_i)$ is a trivial
bundle. Consider for each $i$ a trivialization
$\phi_i:e_i\rightarrow U_i\times V$, the map
${\phi^{-1}}^*(\sigma_{D_{e_i}})$ is an automorphism of the bundle
${\phi_i^{-1}}^*({\pi_S}^*(e))=S_e(U_i)$. This enables to extend
${\phi^{-1}}^*(\sigma_{D_{e_i}})$ to a morphism of the restriction
of the pull-back of $f_n$ to $SU_i$ completing by $1$ on the
diagonal since the orthogonal of $e_U$ is chosen canonically. It
results a morphism $\sigma_D'$ of $C(N)\otimes {\C}^n$, that is as
an element of $C(S^*N)\otimes Gl(V)$. This element defines the
requested element $[\sigma_D']$ of $K_1(C(S^*N)\otimes M_n)\simeq
K_1(C(S^*N))$.

\medskip

Let $B^*N$ be the compactification of $T^*N$ whose fibers are
isomorphic to the unit ball. We consider the bundle $B^*N/T^*N$
that we can identify to the sphere bundle $S^*N$. We have the
following exact sequence
$$0\rightarrow C_0(T^*N)\rightarrow
C(B^*N)\rightarrow C(S^*N)\rightarrow 0.
$$
This sequence gives rise to the following exact sequence in
$K$-theory:
$$
K_1(C(S^*N)\otimes M_n)\rightarrow K_0(C(T^*N)\otimes
M_n)\rightarrow K_0(C(B^*N)\otimes M_n)\rightarrow
K_0(C(S^*N)\otimes M_n)\rightarrow 0.
$$
We can define the boundary $\delta([\sigma_D'])$ which is an
element of $K_0(T^*N\otimes M_n)\simeq K_0(T^*N)$.

\medskip

{\bf Proposition 5.}

{\it The index of $D$ depends only of the class of
$\delta([\sigma_D'])$ in $K_0(T^*N)$.}

\medskip

{\bf Proof.}

We remark that   the image of a symbol by the map $K_1(S^*N\otimes
M_n)\rightarrow K_0(T^*N\otimes M_n)$ is zero, if it is the
restriction of a map defined on $B^*N$. This implies that this
symbol is homotopic to a function which  depends  only of $N$.
Thus the operator it defines is homotopic to a multiplication by a
function whose index is zero.

\bigskip

We end this part by answering the following question: let $N$ be a
manifold endowed with a gerbe $C$, which is the union of two open
sets $U_1$ and $U_2$ such that there exists objects $e_1$ and
$e_2$ of the respective categories $C(U_1)$ and $C(U_2)$. Given
operators $D_{e_1}$ and $D_{e_2}$ on $e_1$ and $e_2$, is it
possible to associate to these operators an element of the
$K$-theory? We do not request any compatibility between $D_{e_1}$
and $D_{e_2}$.

\medskip

 The
vectors bundles $e_1$ and $e_2$ are sub-bundles of the respective
trivial bundles $U_1\times {\C}^n$ and $U_2\times {\C}^p$. Let
$SU_1$ and $SU_2$ be the restriction of the sphere bundle of the
cotangent space of $N$ to $U_1$ and $U_2$, we can canonically
extends the symbols of $D_{e_1}$ and $D_{e_2}$, to respective
automorphisms $\sigma_{D_1}$ and $\sigma_{D_2}$ of $F_1=SU_1\times
{\C}^k$, and $F_2=SU_2\times {\C}^p$, thus elements of
$K_1(C(U_1)\otimes M_n)=K_1(C(U_1))$, and $K_1(C(U_2)\otimes
M_n)=K_1(C(U_2))$, where $C(U_1)$ and $C(U_2)$ are respectively
the set of differentiable functions of $U_1$ and $U_2$.

 The Mayer-Vietoris sequence for algebraic $K$-theory gives rise
 to the sequence
$$
K_2(C(S(U_1\cap U_2))\rightarrow K_1(C(SN))\rightarrow
K_1(C(SU_1))\oplus K_1(C(SU_2))\rightarrow K_1(C((S(U_1\cap U_2)))
$$
if the image of $[\sigma_{D_{e_1}}]+[\sigma_{D_{e_2}}]$ in the
previous sequence is zero, then there exists a class $[\sigma_D]$
of $K_1(SN)$, whose image by the map of the previous exact
sequence is the element $[\sigma_{D_{e_1}}]+[\sigma_{D_{e_2}}]$ of
$K_1(C(SU_1)\oplus K_1(C(SU_2))$. The class $[\sigma_D]$ is not
necessarily unique.

\medskip
\bigskip

{\bf 3. The index formula for operator on gerbes.}

\medskip

Now, we will deduce an index  type theorem for riemannian gerbe.
We know that the Chern character of the cotangent bundle induces
an isomorphism:
$$
K_0(T^*N)\otimes{\R}\rightarrow H_c^{even}(N,{\R})
$$
$$
x\otimes t\rightarrow tch(x)
$$
Consider $Vect(Ind)$ the subspace of $K_0(T^*N)$ generated by
$\sigma_P$, where $P$ is an operator on the riemannian gerbe. It
can be considered as a subspace of $H_c^{even}(N,{\R})$. The map
$Vect(\sigma_P)\rightarrow {\R}$ determined by
$ch([\sigma_P])\rightarrow ind(P)$ can be extended to a linear map
${H_c}^{even}(N,{\R})\rightarrow {\R}$.

The Poincare duality implies the existence of a class $t(N)$ such
that
$$
Ind(P)=\int_{T^*N} ch([\sigma_P])\wedge t(N)
$$

\bigskip

{\bf 4. Applications.}

\bigskip

We will apply now this theory to the problem which has motivated
is construction.

Let $N$ be a riemannian manifold, consider the  Clifford gerbe on
$N$, Let $(U_i)_{i\in I}$ be an open covering of $N$,  the
riemannian connection $w_{C}$ is defined by a family of $so(n)$
$1$-forms $w_i$ on $U_i$ which satisfy
$w_j=ad({g_{ij}}^{-1})w_i+{g_{ij}}^{-1}dg_{ij}$. The covariant
derivative of the Levi-Civita connection is $d+w_i$. We will fixes
an orthogonal basis $(e_1,..,e_n)$ of the tangent space $TU_i$ of
$U_i$. and write $w_i=\sum_{k=1}^{k=n}w_{ik}e_k$.
 we can define the spinorial covariant derivative by setting
$\phi_{ij}=-{1\over 4}w_{ij}$. In the orthogonal basis
$(e_1,..,e_n)$ we have
$\phi_{e_j}=\sum_{k,l}{\phi_{kl}}_{kl}e_ke_l$, with
${\phi}_{kl}=-{\phi}_{lk}$.

\bigskip

{\bf The Dirac operator.}

\bigskip

Let $e_U$ be an object of $Cl(U)$, $U_i$ a trivialization of
$Cl(U)$, we will define
$D_{e_U}=\sum_{k=1}^{k=n}e_i{\nabla_{spin}}_{e_i}$, On each object
$e_U$ of $C(U)$, we have the Lichnerowicz-Weitzenbock formula:
$D^2=\nabla^*\nabla+{1\over 4}s$, where $s$ is the scalar
curvature. In this formula $\nabla^*\nabla$ is the connection
laplacian.
 We will say the global spinor is harmonic if
$D_{e_i}(s_i)=0$, for each $s_i$.

\bigskip

{\bf Proposition 1.}

{\it Suppose that the scalar curvature $s$ of $N$  is strictly
positive, and $N$ is compact then every harmonic global spinor is
$0$.}

\medskip

{\bf Proof.}

Let $\psi$ be an harmonic global spinor, we can represent $\psi$
by a family of spinors $s_i$ defined on an open cover $(U_i)_{i\in
I}$ of $N$, we have on each $U_i$, $D_{e_i}(s_i)=0$, this implies
that ${D_{e_i}}^2(s_i)=0$, we can write
$$
\int_{U_i}<\nabla\nabla^*s_i,s_i>+<{1\over 4}ss_i,s_i>=0
$$
this implies that $\int_{N}s=0$, which  contradicts the fact that
the scalar curvature is strictly positive.

\medskip

\medskip

{\bf Corollary 2.} {\it Suppose that  the sectional curvature of a
compact riemannian manifold is strictly  positive  then the class
$\tau(N)$ associated to the index formula for operators on the
gerbe $Cl(N)$  is zero.}

\bigskip

{\bf Bibliography.}

\bigskip

1. Atiyah, M. Singer, I. The index of elliptic operators, I. Ann.
of Math (87) 1968, 484-530.

\medskip

2. Brylinski, J. Mc Laughin The geometry of degree-four
characteristic classes and line bundles on loop spaces I. Duke
Math. Journal (75) 603-637

\medskip

3. Charles, J. Cours de maitrise de mathematique

\medskip

4. Gilkey, Invariance theory the heat equation, and the
Atiyah-Singer theorem. Perish Publication.

\medskip

5. { A.-O. Kuku}, { Ranks of $K_n$ and $G_n$ of orders and groups
rings of finite groups over integer in number fields }, Journal of
Pure and Applied Algebra., 138 (1999), 39-44.

\medskip
6. Lawson, B. Michelson, L. Spin geometry. Princeton University
Press.

\medskip

7. Murray, M. Singer, M. Gerbes, Clifford modules and the index
theorem. math.DG/0302096

\medskip

8. {M. Nguiffo Boyom},  { Algebres a  associateur symetrique et
algebres reductives}, These Universite de Grenoble.,  (1968).

\medskip

 9.  Van Erp, E. The Atiyah-Singer index theorem,
 $C^*$ algebraic $K-$theory and  quantization
 Master thesis, University of Amsterdam

\end{document}